\documentclass[11pt]{article} 
\begin{document} 
\newtheorem{prop}{Proposition}[section]
\newtheorem{Def}{Definition}[section] \newtheorem{theorem}{Theorem}[section]
\newtheorem{lemma}{Lemma}[section] \newtheorem{Cor}{Corollary}[section]

\title{\bf Modified low regularity well-posedness for the one-dimensional 
Dirac-Klein-Gordon system}
\author{{\bf Hartmut Pecher}\\
Fachbereich Mathematik und Naturwissenschaften\\
Bergische Universit\"at Wuppertal\\
Gau{\ss}str.  20\\
D-42097 Wuppertal\\
Germany\\
e-mail {\tt Hartmut.Pecher@math.uni-wuppertal.de}}
\date{}
\maketitle

\begin{abstract}
The 1D Cauchy problem for the Dirac-Klein-Gordon system is shown to be locally 
well-posed for low regularity Dirac data in $\widehat{H^{s,p}}$ and wave data in 
$\widehat{H^{r,p}} \times \widehat{H^{r-1,p}}$ for $1<p\le 2$ under certain 
assumptions on the parameters r and s, where $\|f\|_{\widehat{H^{s,p}}} := \| 
\langle \xi \rangle^s \widehat{f}\|_{L^{p'}}$ , generalizing the results for 
$p=2$ by Selberg and Tesfahun. Especially we are able to improve the results 
from the scaling point of view with respect to the Dirac part.
\end{abstract}

\renewcommand{\thefootnote}{\fnsymbol{footnote}}
\footnotetext{\hspace{-1.8em}{\it 2000 Mathematics Subject Classification:} 
35Q40, 35L70 \\
{\it Key words and phrases:} Dirac -- Klein -- Gordon system,  
well-posedness, Fourier restriction norm method}
\normalsize 
\setcounter{section}{-1}
\section{Introduction}
Consider the Cauchy problem for the Dirac -- Klein -- Gordon 
equations in one space dimension
\begin{eqnarray}
\label{0.1}
-i\beta\frac{\partial}{\partial t} \psi + i \alpha \beta 
\frac{\partial}{\partial 
x} \psi + M \psi & = & g \phi \psi \\
\label{0.2}
\frac{\partial^2}{\partial t^2} \phi -  \frac{\partial^2}{\partial x^2} \phi + 
m^2 \phi & = & \langle \beta \psi,\psi \rangle_{{\bf C}^2}
\end{eqnarray}
with initial data
\begin{equation}
\psi(x,0)  =  \psi_0(x) \,  , \, \phi(x,0)  =  \phi_0(x) \, , \, \frac{\partial 
\phi}{\partial t}(x,0) = \phi_1(x) \, .
\label{0.3}
\end{equation}
Here $\psi$ is a two-spinor field, i.e. $\psi$ has values in ${\bf C}^2$, and 
$\phi$ is a real-valued function. $\alpha$ and $\beta$ are hermitian 
$(2 \times 2)$ -matrices, which fulfill $\alpha^2 = \beta^2 = I$ , $ \alpha 
\beta + \beta 
\alpha = 0$, e.g. we can choose $\alpha = {0\;\,1 \choose 1\;\,0}$ , $ \beta = 
{1\;\,0\choose0\,-1}$. $M,m$ and $g$ are real constants.

We are interested in local low regularity solutions. This problem was considered 
for data in $L^2$-based Sobolev spaces first by Chadam and Glassey 
\cite{C},\cite{CG} who proved global 
well-posedness for data $\psi_0 \in H^1$ , $\phi_0 \in H^1$ , $\phi_1 \in  L^2$. 
This result was improved by Bournaveas \cite{B} (cf. also Fang \cite{F}) who 
showed the same results for data $\psi_0 \in L^2$ , $\phi_0 \in H^1,$ $\phi_1 
\in L^2$. Local well-posedness was shown by Fang \cite{F1} for data 
$\psi_0 \in H^{-\frac{1}{4}+\epsilon}, $ $ \phi_0 \in H^{\frac{1}{2} +\delta} $ 
, $ \phi_1 \in H^{-\frac{1}{2} + \delta} $ and $0 < \epsilon \le \frac{1}{4},$  
$ 0 < \delta \le 2\epsilon$. Bournaveas and Gibbeson \cite{BG} proved global 
well-posedness for $\psi_0 \in L^2$ , $ \phi_0 \in H^k$ , $ \phi_1 \in H^{k-1}$ 
for $\frac{1}{4} \le k < \frac{1}{2}$. The best known local well-posedness 
result for $L^2$-based Sobolev spaces, namely $\psi_0 \in H^s$, $\phi_0 \in 
H^r,$  $\phi_1 \in H^{r-1}$ was given by Selberg and Tesfahun \cite{ST} who 
assumed $s > -1/4$, $ r > 0$ , $|s| \le r \le 1+s$ . This result was also shown 
to be optimal within the used method, namely the Bougain - Klainerman - Machedon 
Fourier restriction norm method. This also improved earlier results of Machihara 
\cite{M} and the author \cite{P}. The most recent results of \cite{ST} and 
\cite{P} were obtained using the null structure of $\langle \beta \psi,\psi 
\rangle$ of the wave part, which is also hidden (by a duality argument) in the 
Dirac part of the system. This fact was first detected by d'Ancona, Foschi and 
Selberg \cite{AFS}, who showed local well-posedness for data $\psi_0 \in 
H^{\epsilon}$ , $\phi_0 \in H^{\frac{1}{2}+\epsilon}$ , $ \phi_1 \in 
H^{-\frac{1}{2}+\epsilon}$ (for $\epsilon > 0$) in the (3+1)-dimensional case, 
being arbitrarily close to the minimal regularity predicted by scaling 
($\epsilon = 0$). In contrast to this (3+1)-dimensional result the best result 
in one space dimension by \cite{ST} mentioned above is far away from the 
regularity predicted by scaling, namely $\psi_0 \in H^{-1}$ , $\phi_0 \in 
H^{-1/2}$ , $ \phi_1 \in H^{-3/2}$ .

The aim of the present paper is to close this gap as far as possible. In order 
to achieve this we leave the $H^s$-scale of the data spaces. This was suggested 
for nonlinear Schr\"odinger equations by Cazenave, Varga, and Vilela \cite{CVV} 
and Vargas and Vega \cite{VV}. This method, a modified Fourier restriction norm 
method, was systematically introduced by Gr\"unrock in \cite{G1}, where he 
applied it to the modified KdV-equation and in \cite{G2}, where he was able to 
show local well-posedness in the case of the cubic nonlinear Schr\"odinger 
equation
$ iu_t + u_{xx} + |u|^2 u = 0 $ , $ u(0) = u_0 $ , for data $u_0 \in 
\widehat{H^{s,p}}$, where
\begin{equation}
\label{0.0}
\|u_0
\|_{{\widehat{H^{s,p}}}} := \| \langle \xi \rangle^s \widehat{u_0} 
\|_{H^{p'}_{\xi}} 
\end{equation} 
$ 1/p + 1/p' = 1 $ , if $ s \ge 0$ and $1<p<\infty$ , and global well-posedness 
for $2 \ge p \ge 5/3,$  $u_0 \in \widehat{H^{0,p}}$ .

In the present paper we prove local well-posedness for the Dirac-Klein-Gordon 
system with data $\psi_0 \in \widehat{H^{s,p}}$ , $\phi_0 \in \widehat{H^{r,p}}$ 
, $\phi_1 \in \widehat{H^{r-1,p}}$ for $1<p \le 2$ under suitable assumptions on 
s and r. These results allow to improve the $L^2$-based results from the scaling 
point of view. More precisely we need the following conditions:
$$ s > -\frac{1}{2} + \frac{1}{2p} \, , \, r \ge |s| \, , \, r > \frac{2}{p} - 1 
\, , \, r \le 1+s \, , $$
which reduce to the results of \cite{ST} in the case $p=2$ . Especially we are 
able to choose $p=1+$ , $s=0$ , $r=1$ , leading to $\psi_0 \in 
\widehat{H^{0,1+}}$ , which scales like $H^{-\frac{1}{2}+,2}$ and to 
$(\phi_0,\phi_1) \in \widehat{H^{1,1+}} \times \widehat{H^{0,1+}}$, which scales 
like $H^{\frac{1}{2}+,2} \times H^{-\frac{1}{2}+,2}$ . Thus the result by 
\cite{ST},\cite{P} in the case $p=2$, namely $\psi_0 \in 
H^{-\frac{1}{4}+\epsilon,2}$ , $ (\phi_0,\phi_1) \in H^{\frac{1}{4}-\epsilon,2} 
\times H^{-\frac{3}{4}-\epsilon,2}$ , $\frac{1}{4} > \epsilon > 0$ , is improved 
for the Dirac part at the expense of weakening the result for the wave part.

This paper is organized as follows: We diagonalize the system like \cite{ST} by 
using the projections $P_{\pm}$ onto the eigenspaces of $-i \alpha 
\frac{\partial}{\partial x}$ and splitting $\psi$ as a sum $P_+ \psi + P_- 
\psi$. We also split $\phi$ as a sum $\phi = \phi_+ + \phi_-$, where the half 
waves $\phi_+$ and $\phi_-$ are defined in the usual way. Then we analyze the 
components of the nonlinearity $\langle \beta \psi,\psi ' \rangle$, namely 
$\langle \beta P_{\pm} \psi, P_{\pm} \psi ' \rangle$ for all combinations of 
signs. It turns out that $\langle \beta P_+ \psi, P_- \psi ' \rangle$ and 
$\langle \beta P_- \psi, P_+ \psi ' \rangle$ vanish. Then we examine which 
bilinear estimates for the nonlinear terms are necessary for local 
well-posedness in the framework of the $X^{l,b}_p$ - spaces (for a definition 
cf. (\ref{0.a}),(\ref{0.b})). These are proven in Proposition 
\ref{Proposition1.1}.
The results are summarized in Theorem \ref{Theorem} and the Remark to Theorem 
\ref{Theorem}. 

We recall the modified Fourier restriction norm method in the following. For 
details we refer to the paper of A. Gr\"unrock (cf. \cite{G1}, Chapter 2).
Our solution spaces are the Banach spaces
\begin{equation}
\label{0.a}
 X_p^{l,b} := \{ f \in {\cal S}'({\bf R}^2) : \|f\|_{X^{l,b}_p} < \infty\} \, 
, 
\end{equation}
where $l,b \in {\bf R}$ , $1<p<\infty$ , $1/p + 1/p' = 1$ and
\begin{equation}
\label{0.b}
\|f\|_{X^{l,b}_p} := \left( \int d\xi d\tau \langle \xi \rangle^{lp'} 
\langle 
\tau + \phi(\xi) \rangle^{bp'} |\hat{f}(\xi,\tau)|^{p'} \right)^{1/p'} \, , 
\end{equation}
where $\phi: {\bf R} \to {\bf R} $ is a given smooth function of polynomial 
growth. We denote by $\tilde{f}$ or ${\cal F}f$ the Fourier transform with 
respect to space and time. The dual space of $X^{l,b}_p$ is $X^{-l,-b}_{p'}$ , 
and the Schwartz 
space 
is dense in $X^{l,b}_p$ . We have 
$(X^{l_0,b_0}_{p_0},X^{l_1,b_1}_{p_1})_{[\Theta]} = X^{l,b}_p $ , where 
$[\Theta]$ denotes the complex interpolation method, and for $l_0,l_1,b_0,b_1 
\in {\bf R}$ , $ 1 < p_0,p_1 \le \infty $ , $ \Theta \in [0,1] $ we have 
$l=(1-\Theta)l_0+\Theta l_1 $ , $ b = (1-\Theta)b_0 + \Theta b_1 $ . $ 
\frac{1}{p} = \frac{1-\Theta}{p_0} + \frac{\Theta}{p_1}$ . The embedding 
$X^{l,b}_p \subset C^0({\bf 
R},\widehat{H^{l,p}})$ is true for $b>1/p$ (recall the definition (\ref{0.0})). 
We have
$$ \|f\|_{X^{l,b}_p} = \left(\int d\xi d\tau \langle \xi \rangle^{lp'} \langle 
\tau \rangle^{bp'} \left| {\cal F}(U_{\phi}(- \cdot) f)(\xi,\tau) 
\right|^{p'} \right)^{1/p'} \, .$$ Here $U_{\phi}(t) := 
e^{-it\phi(-i\partial_x)}$.
For any $\psi \in C^{\infty}_0({\bf R}_t)$ one has
$$ \| \psi U_{\phi}(t) u_0\|_{X^{l,b}_p} \le c_{\psi} 
\|u_0\|_{\widehat{H^{l,p}}} \, .$$
If $v$ is a solution of the inhomogeneous problem
$$ i v_t - \phi(-i\partial_x)v = F \, , \, v(0) = 0 $$
and $\psi \in C^{\infty}_0({\bf R}_t)$ with $supp \, \psi \subset (-2,2) $ , $ 
\psi \equiv 1 $ on $[-1,1]$ , $ \psi(t) = \psi(-t) $ , $ \psi(t) \ge 0 $ , $ 
\psi_{\delta}(t) := \psi(\frac{t}{\delta}) $ , $ 0<\delta \le 1$ , we have for 
$1<p<\infty$ , $ b'+1 \ge b \ge 0 \ge b' > - 1/p' $ :
$$ \|\psi_{\delta} v\|_{X^{l,b}_p} \le c \delta^{1+b'-b} \|F\|_{X^{l,b'}_p} \, 
. 
$$
We also use the localized 
spaces
$$ X^{l,b}_p[0,T] := \{ f = \bar{f}_{|[0,T] \times {\bf R}} : \bar{f} \in 
X^{l,b}_p \} \, , $$
where
$$
\|f\|_{X^{l,b}_p[0,T]} := \inf \{ \|\bar{f}\|_{X^{l,b}_p} : f 
= \bar{f}_{|[0,T]\times{\bf R}}\} \, . $$
For the case $\phi(\xi) = \pm \xi$ we use the notation $X^{l,b}_{\pm p}$ and 
$\|\psi\|_{X^{l,b}_{\pm p}} = \| \langle \xi \rangle^l \langle \tau \pm \xi 
\rangle^b \tilde{\psi}(\xi,\tau)\|_{L^{p'}_{\xi \tau}} $ , whereas for the case 
$\phi(\xi) = \pm |\xi|$ we denote the space by $Y^{l,b}_{\pm p}$ and 
$\|\psi\|_{Y^{l,b}_{\pm p}} = \| \langle \xi \rangle^l \langle \tau \pm |\xi| 
\rangle^b \tilde{\psi}(\xi,\tau)\|_{L^{p'}_{\xi \tau}} $.

Especially we use \cite{G1}, Theorem 2.3, which we repeat for convenience. 
\begin{theorem}
\label{Theorem0.1}
Consider the Cauchy problem
\begin{equation}
\label{0*}
i u_t -  \phi(-i\partial_x) u = N(u) \quad , \quad u(0) = u_0 \in 
\widehat{H^{s,p}} \, , 
\end{equation}
where $N$ is a nonlinear function of $u$ and its spatial derivatives. Assume 
for 
given $s \in {\bf R} $ , $1<p<\infty$ , $\alpha \ge 1$ there exist $ b > 1/p $ 
, 
$b-1 < b' \le 0$ such that the estimates
$$ \|N(u)\|_{X^{s,b'}_p} \le c \|u\|_{X^{s,b}_p}^{\alpha} $$
and
$$ \|N(u)-N(v)\|_{X^{s,b'}_p} \le c (\|u\|_{X^{s,b}_p}^{\alpha -1} + 
\|v\|_{X^{s,b}_p}^{\alpha -1}) \|u-v\|_{X^{s,b}_p}$$
are valid. Then there exist $T = T(\|u_0\|_{\widehat{H^{s,p}}}) > 0 $ and a 
unique solution $u \in X^{s,b}_p[0,T]$ of (\ref{0*}). This solution belongs to 
$C^0([0,T],\widehat{H^{s,p}})$ , and the mapping $u_0 \mapsto u$ , 
$\widehat{H^{s,p}} \to X^{s,b}_p(0,T_0)$ is locally Lipschitz continuous for 
any 
$T_0 < T$ .
\end{theorem}
We use the notation $\langle \lambda \rangle := (1+\lambda^2)^{1/2}$ , and $a 
\pm$ to denote a number slightly larger (resp., smaller) than $a$ .

\section{Local well-posedness}
First we transform our system (\ref{0.1}),(\ref{0.2}) into a first order system 
(in 
t) in diagonal form.

Multiplying the Dirac equations by $\beta$ leads to
\begin{eqnarray*}
-i\frac{\partial}{\partial t} \psi - i \alpha  \frac{\partial}{\partial x} \psi 
+ 
M \beta \psi & = & g \phi \beta \psi \\
\frac{\partial^2}{\partial t^2} \phi -  \frac{\partial^2}{\partial x^2} \phi + 
m^2 \phi & = & \langle \beta \psi,\psi \rangle_{{\bf C}^2} \, .
\end{eqnarray*}
Following the paper of Selberg and Tesfahun we diagonalize the system by 
defining the projections
$$ P_{\pm} := \frac{1}{2} {\;\, 1 \;\,\pm 1 \choose \pm 1\;\,1} \, . $$ Then we 
have $\psi = \psi_+ + \psi_-$ 
with $\psi_{\pm} := P_{\pm}\psi $. Using the identities
$ \alpha = P_+ - P_- ,$ $ P_{\pm}^2 = P_{\pm} $ and $ P_{\pm}P_{\mp} = 0$ we get 
by application of $P_{\pm}$ to the Dirac equation 
\begin{eqnarray*}
P_{\pm}(-i\frac{\partial}{\partial t}\psi - i\alpha\frac{\partial}{\partial 
x}\psi) & = & -i\frac{\partial}{\partial t}P_{\pm}\psi -i 
\frac{\partial}{\partial x} P_{\pm}(P_+ - P_-)\psi \\
& = & -i\frac{\partial}{\partial t} \psi_{\pm} \mp i \frac{\partial}{\partial x} 
\psi_{\pm}
\end{eqnarray*}
and thus the Dirac equations are transformed into
\begin{eqnarray}
\nonumber 
(-i\frac{\partial}{\partial t} \mp i \frac{\partial}{\partial x})\psi_{\pm} & = 
& -M\beta 
P_{\mp}(\psi_+ + \psi_-) + g P_{\pm}(\phi \beta \psi) \\
\label{1.1}
& = & -M\beta \psi_{\mp} + g P_{\pm}(\phi \beta (\psi_+ + \psi_-)) \, ,
\end{eqnarray}
where we also used $P_{\pm} \beta = \beta P_{\mp}$ .
We also split the function $\phi$ into the sum $\phi = \frac{1}{2}(\phi_+ + 
\phi_-)$ , where
$$ \phi_{\pm} := \phi \pm iA^{-\frac{1}{2}} \frac{\partial \phi}{\partial t} 
\quad , \quad A:= - \frac{\partial^2}{\partial x^2} + m^2 \, $$
Here we assume $m > 0$ and in fact $m=1$. Otherwise we artificially add a term 
$(1-m^2)\phi$ on both sides of the equation at the expense of having an 
additional linear term $c_0 \phi$ in the inhomogeneous part which can easily be 
taken care of. We easily calculate
\begin{equation}
\label{1.2}
(i \frac{\partial}{\partial t} \mp A^{\frac{1}{2}})\phi_{\pm}  =  
 \mp A^{-\frac{1}{2}}( \langle \beta \psi,\psi\rangle_{{\bf C}^2} + c_0 
\phi) \, . 
\end{equation}
The initial conditions are transformed into
\begin{equation}
\label{1.3}
\psi_{\pm}(0,x) = P_{\pm}\psi_0(x) \quad , \quad \phi_{\pm}(0,x) = 
\phi_0(x) 
\pm iA^{-\frac{1}{2}} \phi_1(x) \, . 
\end{equation}

The following (slightly modified) system of integral equations  
belongs to our Cauchy problem (\ref{1.1}),(\ref{1.2}),(\ref{1.3}), where 
$U_{\pm}(t)$ denotes the evolution operator of the equation 
$(\frac{\partial}{\partial t} \pm \frac{\partial}{\partial x})u = 0$.
\begin{eqnarray}
\label{a}
\psi_{\pm}(t) & = &  U_{\pm}(t)\psi_{\pm}(0) \\
\nonumber
& & -i g \int_0^t U_{\pm}(t-s) P_{\pm} (\frac{1}{2} 
(\phi_+(s)+\phi_-(s)) 
\beta(P_+ \psi_+(s) \\
\nonumber
& & \quad \quad + P_- \psi_-(s))) ds + iM \int_0^t U_{\pm}(t-s) \beta 
P_{\mp}\psi_{\mp}(s) ds \\
\label{b}
\phi_{\pm}(t) & = & e^{\mp itA^{\frac{1}{2}}} \phi_{\pm}(0)\\
\nonumber
& & \pm i \int_0^t e^{\mp i(t-s)A^{\frac{1}{2}}} A^{-\frac{1}{2}} \langle 
\beta(P_+ 
\psi_+(s) + P_- \psi_-(s)),P_+ \psi_+(s) \\
\nonumber
& & \quad \quad + P_- \psi_-(s) \rangle ds  \pm i c_0 \int_0^t 
e^{\mp i(t-s)A^{\frac{1}{2}}} A^{-\frac{1}{2}} (\phi_+(s) + \phi_-(s)) ds
\end{eqnarray}
We remark that any solution of this system automatically fulfills 
$P_{\pm}\psi_{\pm} = \psi_{\pm},$ because applying $P_{\pm}$ to the 
right hand side of the equations for $\psi_{\pm}$ gives 
$P_{\pm}\psi_{\pm}(0) = P_{\pm}P_{\pm}\psi_0 = P_{\pm}\psi_0 
= \psi_{\pm}(0)$ , and the integral terms also remain unchanged, because 
$P_{\pm}^2 = P_{\pm}$ and $P_{\pm}\beta P_{\mp}\psi_{\mp}(s) = \beta 
P_{\mp}\psi_{\mp}(s)$ . Thus $P_{\pm}\psi_{\pm}$ 
can be replaced by $\psi_{\pm}$ on the right hand sides, thus the system of 
integral equations reduces exactly to the one belonging to our Cauchy problem 
(\ref{1.1}),(\ref{1.2}),(\ref{1.3}).
 
Let now data be given with 
$$ \psi_0 \in \widehat{H^{s,p}} \, , \, \phi_0 \in \widehat{H^{r,p}} \, , \, 
\phi_1 \in \widehat{H^{r-1,p}} \, . $$
This implies $ \psi_{\pm}(0) \in \widehat{H^{s,p}}$ and $\phi_{\pm}(0) \in 
\widehat{H^{r,p}}$ . In order to construct a solution of the integral equations 
for 
$t \in [0,T]$ with $ \psi_{\pm} \in X_{\pm p}^{s,\sigma}[0,T]$ and $\phi_{\pm} 
\in 
Y_{\pm p}^{r,\rho}[0,T],$ where $1 > \sigma,\rho > 1/p$ ,  we may apply Theorem 
\ref{Theorem0.1}, because its generalization from the case of a single equation 
to a system is evident. Thus we only have to 
show the following estimates for the nonlinearities.

Concerning (\ref{a}) we need
\begin{equation}
\label{****}
\| P_{\pm}(\phi \beta P_{[\pm]} \psi) 
\|_{X^{s,\sigma -1+\epsilon}_{\pm p}} \le c \| \phi 
\|_{Y^{r,\rho}_{+ p}} \| \psi 
\|_{X^{s,\sigma}_{[\pm] p}} 
\end{equation}
and the same estimates with $Y^{r,\rho}_{+ p}$ 
replaced by $Y^{r,\rho}_{- p}$. 
$[\pm]$ denotes a sign independent of $\pm$.
By duality this is equivalent to
$$ \left| \int \int \langle P_{\pm}(\phi \beta P_{[\pm]} \psi),\psi' 
\rangle dx dt \right|  \le c \|\phi\|_{Y^{r,\rho}_{+ p}} 
\|\psi\|_{X^{s,\sigma}_{[\pm] p}} 
\|\psi'\|_{X^{-s,1-\sigma -\epsilon}_{\pm p'}} \, . $$
The left hand side equals
$$ \left| \int \int \phi \langle \beta P_{[\pm]} \psi,P_{\pm} \psi' 
\rangle dx dt \right| \, , $$
which can be estimated by
$$  \|\phi\|_{Y^{r,\rho}_{+ p}} \| \langle \beta P_{[\pm]}(D) 
\psi , P_{\pm} \psi' \rangle \|_{Y^{-r,-\rho}_{+ p'}} \, . $$
Thus (\ref{****}) is fulfilled if
\begin{equation}
\label{***}
\| \langle \beta P_{[\pm]} \psi , P_{\pm} \psi' \rangle 
\|_{Y^{-r,-\rho}_{+ p'}} \le c 
\|\psi\|_{X^{s,\sigma}_{[\pm] p}} 
\|\psi'\|_{X^{-s,1-\sigma - \epsilon}_{\pm p'}}
\end{equation}
Concerning (\ref{b}) we have to show
\begin{equation}
\label{**}
\| \langle \beta P_{[\pm]} \psi, P_{\pm} \psi' \rangle 
\|_{Y^{r-1,\rho -1+\epsilon}_{+ p}} \le c 
\|\psi\|_{X^{s,\sigma}_{[\pm] p}} 
\|\psi'\|_{X^{s,\sigma}_{\pm p}} \, .
\end{equation}
We also need the same estimates with $Y_+$ replaced by $Y_-$ .

The linear terms in the integral equations can easily be treated as follows:
Let $\varphi$ be a $t$-dependent $C^{\infty}$ - function with $ \varphi = 1 $ on 
$[0,T]$ and $supp \, \varphi \subset [0,2T]$. Then, with $J^s$ being the 
multiplier with symbol $\langle \xi \rangle^s$ , and using the embedding 
(\ref{2.2}) we have:
\begin{eqnarray*}
\|\psi_{\pm} \|_{X^{s,\sigma -1+\epsilon}_{\pm p}[0,T]} & \le & \|\varphi 
\overline{\psi_{\pm}}\|_{X^{s,\sigma -1+\epsilon}_{\pm p}} \le \|\varphi 
\overline{\psi_{\pm}}\|_{X^{s,0}_{\pm p}} \\ 
& = & \|\varphi J^s \overline{\psi_{\pm}} \|_{\widehat{L^p_{xt}}} \le 
\|\varphi\|_{\widehat{L^p_t}} \|J^s \overline{\psi_{\pm}} 
\|_{\widehat{L^p_x}(\widehat{L^{\infty}_t})} \\ 
& \le & \|\varphi\|_{\widehat{L^p_t}} \|J^s \overline{\psi_{\pm}} 
\|_{X^{0,\sigma}_{\pm p}} \le c T^{1/p} 
\|\overline{\psi_{\pm}}\|_{X^{s,\sigma}_{\pm p}}
\end{eqnarray*}
for any $\overline{\psi_{\pm}}$ with $\overline{\psi_{\pm}}_{|[0,T]} = 
\psi_{\pm}$ .
Here we used Hausdorff-Young to estimate $\|\varphi\|_{\widehat{L^p_t}}$ by 
$cT^{1/p}$. Thus:
$$ \|\psi_{\pm}\|_{X^{s,\sigma -1+\epsilon}_{\pm p}[0,T]} \le c T^{1/p} 
\|\psi_{\pm}\|_{X^{s,\sigma}_{\pm p}[0,T]}$$
and similarly
$$ \|A^{-1/2} \phi_{\pm}\|_{Y^{0,\rho -1+\epsilon}_{\pm p}[0,T]} \le c T^{1/p} 
\|\phi_{\pm}\|_{Y^{r,\rho}_{\pm p}[0,T]} \, . $$
The bilinear form has a null structure, we namely have
$$ \langle \beta P_{\pm} \psi,P_{\pm} \psi' \rangle = \langle P_{\mp} \beta 
\psi,P_{\pm} \psi' \rangle = 0 \, , $$
so that in order to prove (\ref{***}) and (\ref{**}) it remains to show
\begin{eqnarray}
\label{*1}
\| \langle \beta P_+ \psi , P_- \psi' \rangle 
\|_{Y^{-r,-\rho}_{\pm p'}} & \le & c 
\|\psi\|_{X^{s,\sigma}_{+ p}} 
\|\psi'\|_{X^{-s,1-\sigma - \epsilon}_{- p'}} \\
\label{*2}
\| \langle \beta P_- \psi , P_+ \psi' \rangle 
\|_{Y^{-r,-\rho}_{\pm p'}} & \le & c 
\|\psi\|_{X^{s,\sigma}_{- p}} 
\|\psi'\|_{X^{-s,1-\sigma - \epsilon}_{+ p'}} \\
\label{**1}
\| \langle \beta P_+ \psi, P_- \psi' \rangle 
\|_{Y^{r-1,\rho -1+\epsilon}_{\pm p}} & \le & c 
\|\psi\|_{X^{s,\sigma}_{+ p}} 
\|\psi'\|_{X^{s,\sigma}_{- p}} \\
\label{**2}
\| \langle \beta P_- \psi, P_+ \psi' \rangle 
\|_{Y^{r-1,\rho -1+\epsilon}_{\pm p}} & \le & c 
\|\psi\|_{X^{s,\sigma}_{- p}} 
\|\psi'\|_{X^{s,\sigma}_{+ p}} \, .
\end{eqnarray}
In order to prove (\ref{**1}) we have to show
$$ \left| \int \int \langle \beta P_+ \psi , P_- \psi' 
\rangle \overline{\phi} dx dt \right| \le c 
\|\phi\|_{Y^{1-r,1-\rho - \epsilon}_{\pm p'}} 
\|\psi\|_{X^{s,\sigma}_{+ p}} 
\|\psi'\|_{X^{s,\sigma}_{- p}} \, . $$
The left hand side equals
$$ \left| \int \int_*  \langle \beta P_+  
\tilde{\psi}(\xi_1,\tau_1) , P_- \tilde{\psi}'(-\xi_2,-\tau_2) \rangle 
\overline{\tilde{\phi}}(\xi,\tau) d\xi_1 d\xi_2 d\tau_1 d\tau_2 \right| \, , $$
where * denotes the region $\xi_1 + \xi_2 = \xi$ , $ \tau_1 + \tau_2 = \tau $ .

Defining now 
\begin{eqnarray*}
\tilde{v}_1(\xi_1,\tau_1) & := & \langle \xi_1 \rangle^s \langle \tau_1 + 
\xi_1 \rangle^{\sigma}\tilde{\psi}(\xi_1,\tau_1) \\
\tilde{v}_2(\xi_2,\tau_2) & := & \langle \xi_2 \rangle^s \langle \tau_2 - 
\xi_2 \rangle^{\sigma}\tilde{\psi'}(\xi_2,\tau_2) \\
\tilde{\varphi}(\xi,\tau) & := & \langle \xi \rangle^{1-r} \langle \tau \pm 
|\xi| 
\rangle^{1-\rho - \epsilon}\tilde{\phi}(\xi,\tau) 
\end{eqnarray*}
we have
$$ \|\psi\|_{X^{s,\sigma}_{+ p}} = \|v_1\|_{\widehat{L^p_{xt}}} \, , \,
\|\psi'\|_{X^{s,\sigma}_{- p}} = \|v_2\|_{\widehat{L^p_{xt}}} \, , \, 
\|\phi\|_{Y^{1-r,1-\rho - \epsilon}_{\pm p'}} = 
\|\varphi\|_{\widehat{L^{p'}_{xt}}} \, . $$ 
Thus we have to show
\begin{eqnarray*}
\int \int_* \frac{ 
|\tilde{v}_1(\xi_1,\tau_1)\tilde{v}_2(-\xi_2,-\tau_2) 
\overline{\tilde{\varphi}}(\xi,\tau)|}{ \langle \xi_1 \rangle^s \langle 
\tau_1 + \xi_1 \rangle^{\sigma} \langle \xi_2 \rangle^s \langle \tau_2 - \xi_2 
\rangle^{\sigma} \langle \xi \rangle^{1-r} \langle \tau \pm |\xi| 
\rangle^{1-\rho - \epsilon}} d\xi_1 d\xi_2 
d\tau_1 d\tau_2  & & \\ 
 \le  c \|v_1\|_{\widehat{L^p}} \|v_2\|_{\widehat{L^p}} 
\|\varphi\|_{\widehat{L^{p'}}} \, . & & 
\end{eqnarray*}
Define $\sigma_1^+ := \tau_1 + \xi_1$ , $ \sigma_2^- := \tau_2 - \xi_2$ , $ 
\sigma_{\pm} := \tau \pm |\xi|$ . \\
{\bf Case 1:} $|\sigma_{\pm}|$ dominant, i.e. $|\sigma_{\pm}| \ge 
|\sigma_1^+|,|\sigma_2^-|$ . \\
{\bf a:} $|\xi_1| << |\xi_2| $ ( $ \Rightarrow |\xi| \sim |\xi_2| $ ) . (The 
case $|\xi_2|<<|\xi_1|$ can be treated similarly.) \\
In this case we have by the algebraic inequality in Lemma \ref{Lemma2.1} the 
estimate $ \langle \sigma_{\pm} \rangle \ge c \langle \xi_1 \rangle$ , so that 
it remains to estimate (provided $\rho \le 1-\epsilon$):
\begin{eqnarray*}
\int \int_* \frac{ 
|\tilde{v}_1(\xi_1,\tau_1)\tilde{v}_2(-\xi_2,-\tau_2) 
\overline{\tilde{\varphi}}(\xi,\tau)|}{ \langle \xi_1 
\rangle^{s+1-\rho-\epsilon} \langle 
\sigma_1^+ \rangle^{\sigma} \langle \xi_2 \rangle^{s+1-r} \langle \sigma_2^- 
\rangle^{\sigma} } d\xi_1 d\xi_2 
d\tau_1 d\tau_2  & & \\
\le c \int \int_* \frac{ 
|\tilde{v}_1(\xi_1,\tau_1)\tilde{v}_2(-\xi_2,-\tau_2) 
\overline{\tilde{\varphi}}(\xi,\tau)|}{ \langle 
\sigma_1^+ \rangle^{\sigma} \langle \sigma_2^- 
\rangle^{\sigma}} d\xi_1 d\xi_2 
d\tau_1 d\tau_2 \, ,
\end{eqnarray*}
where we made the assumptions
\begin{equation}
\label{1}
s+1-\rho-\epsilon \ge 0
\end{equation}
and
\begin{equation}
\label{2}
s+1-r \ge 0 \, .
\end{equation}
Using Corollary \ref{Corollary2.1} below we get the desired bound.\\
{\bf b:} $|\xi_1| \sim |\xi_2| $ \\
Using the estimate $ \langle \sigma_{\pm} \rangle \ge c \langle \xi_1 \rangle$ 
again we arrive at
$$ \int \int_* \frac{ 
|\tilde{v}_1(\xi_1,\tau_1)\tilde{v}_2(-\xi_2,-\tau_2) 
\overline{\tilde{\varphi}}(\xi,\tau)|}{ \langle \xi_1 
\rangle^{2s+1-\rho-\epsilon} \langle \xi \rangle^{1-r} \langle 
\sigma_1^+ \rangle^{\sigma} \langle \sigma_2^- 
\rangle^{\sigma} } d\xi_1 d\xi_2 
d\tau_1 d\tau_2 \, . $$
If $r>1$ we further have $\langle \xi \rangle^{r-1} \le c(\langle \xi_1 
\rangle^{r-1} + \langle \xi_2 \rangle^{r-1})$ and get the bound
$$ \int \int_* \frac{ 
|\tilde{v}_1(\xi_1,\tau_1)\tilde{v}_2(-\xi_2,-\tau_2) 
\overline{\tilde{\varphi}}(\xi,\tau)|}{ \langle \xi_1 
\rangle^{2s+1-\rho-\epsilon +1-r} \langle 
\sigma_1^+ \rangle^{\sigma} \langle \sigma_2^- 
\rangle^{\sigma} } d\xi_1 d\xi_2 
d\tau_1 d\tau_2 \, . $$
Using again the bilinear estimate of Corollary \ref{Corollary2.1} we get the 
desired bound provided
\begin{equation}
\label{3}
2s+1-\rho-\epsilon \ge 0
\end{equation}
and
\begin{equation}
\label{4}
2s+1-\rho-\epsilon+1-r \ge 0 \, .
\end{equation}
{\bf Case 2:} $|\sigma_1^+|$ dominant and $|\xi_2| \ge |\xi_1|$ . (The other 
cases: $|\sigma_2^-|$ dominant and/or $|\xi_2| \le |\xi_1|$ are similar.) \\
Using the algebraic inequality $|\sigma_1^+| \ge c|\xi_1|$ we have to estimate
$$ \int \int_* \frac{ 
|\tilde{v}_1(\xi_1,\tau_1)\tilde{v}_2(-\xi_2,-\tau_2) 
\overline{\tilde{\varphi}}(\xi,\tau)|}{ \langle \xi_1 \rangle^{s+\sigma} \langle 
\xi_2 \rangle^s \langle \sigma_2^- 
\rangle^{\sigma} \langle \xi \rangle^{1-r} \langle \sigma_{\pm} 
\rangle^{1-\rho-\epsilon} } d\xi_1 d\xi_2 
d\tau_1 d\tau_2 \, . $$
We apply Proposition \ref{Proposition2.1} with $a=s+\sigma$ , $b=s$ , $c=1-r$, 
$\alpha =0$ , $\beta = \sigma$ , $ \gamma = 1-\rho-\epsilon$ and get the desired 
estimate provided
\begin{eqnarray}
\label{5}
2s + \sigma + 1-r & > & \frac{1}{p} + \frac{1}{p} + \frac{1}{p'} - 1 = 
\frac{1}{p} \\
\label{6}
2s + \sigma & \ge & 0 \\
\label{7}
s+\sigma+1-r & \ge &  0 \\
\label{8}
s+1-r & \ge & 0 \\
\label{9}
\sigma +1-\rho-\epsilon & > & \frac{1}{p} + \frac{1}{p} + \frac{1}{p'} - 1 = 
\frac{1}{p} \\
\label{10}
\rho & \le & 1-\epsilon \, .
\end{eqnarray}
We conclude that (\ref{**1}) (and similarly (\ref{**2})) holds, if (\ref{1}) -  
(\ref{10}) are satisfied.

Next we have to prove (\ref{*1}). Similarly as before we have to show
\begin{eqnarray*}
\int \int_* \frac{ 
|\tilde{v}_1(\xi_1,\tau_1)\tilde{v}_2(-\xi_2,-\tau_2) 
\overline{\tilde{\varphi}}(\xi,\tau)|}{ \langle \xi_1 \rangle^s \langle 
\sigma_1^+ \rangle^{\sigma} \langle \xi_2 \rangle^{-s} \langle \sigma_2^- 
\rangle^{1-\sigma-\epsilon} \langle \xi \rangle^r \langle \sigma_{\pm} 
\rangle^{\rho}} d\xi_1 d\xi_2 
d\tau_1 d\tau_2  & & \\ 
 \le  c \|v_1\|_{\widehat{L^p}} \|v_2\|_{\widehat{L^{p'}}} 
\|\varphi\|_{\widehat{L^p}} \, . & & 
\end{eqnarray*}
We use the algebraic inequality of Lemma \ref{Lemma2.1} and remark that
$\min (\sigma,\rho,1-\sigma-\epsilon) = 1-\sigma-\epsilon$ , because 
$\sigma,\rho > 1/p > 1/2$ for $1<p\le 2$. Assuming
\begin{equation}
\label{11}
\sigma \le 1-\epsilon
\end{equation}
we can therefeore replace one of the expressions $\langle \sigma_1^+ 
\rangle^{1-\sigma-\epsilon}$ , $\langle \sigma_2^+-\rangle^{1-\sigma-\epsilon}$ 
or $\langle \sigma_{\pm} \rangle^{1-\sigma-\epsilon}$ by $\langle \xi_1 
\rangle^{1-\sigma-\epsilon}$ or $\langle \xi_2 \rangle^{1-\sigma-\epsilon}$ . In 
any case the conditions of Proposition \ref{Proposition2.1} are satisfied to 
give the desired bound provided the following assumptions are made:
\begin{eqnarray}
\label{12}
s-s+r+1-\sigma-\epsilon & > & \frac{1}{p} \\
\label{13}
s-s+1-\sigma-\epsilon & \ge & 0 \\
\label{14}
s+r & \ge & 0 \\
\label{15}
-s+r & \ge & 0 \\
\label{16}
s+1-\sigma-\epsilon +r & \ge & 0 \\
\label{17}
-s+1-\sigma-\epsilon +r & \ge & 0 \, .
\end{eqnarray}
The sum of the exponents of the remaining $\sigma$-modules is $\rho + \sigma$ 
which has to be larger than $1/p$. This is trivially satisfied. \\
(\ref{*2}) is proven completely analogously.

We summarize our results in the following
\begin{prop}
\label{Proposition1.1}
Let $1 < p \le 2$. The inequalities 
(\ref{*1}),(\ref{*2}),(\ref{**1}),(\ref{**2}) are satisfied with suitable $1 > 
\sigma,\rho > 1/p$, if the following conditions hold:
\begin{eqnarray}
\label{18}
s & > & -\frac{1}{2} + \frac{1}{2p} \\
\label{19}
r & \le & 1+s \\
\label{20}
r & \ge & |s| \\
\label{21}
r & > & \frac{2}{p} - 1 \, .
\end{eqnarray}
\end{prop}
{\bf Proof:} We only have to check the conditions (\ref{1}) - (\ref{17}). Choose 
$ \rho = \frac{1}{p} + \epsilon $ , $\epsilon > 0$ small. Then: \\
(\ref{1}) is satisfied, because by use of (\ref{18}) we have $ s+1-\rho-\epsilon 
> \frac{1}{2}-\frac{1}{2p}-2\epsilon > 0$. \\
(\ref{2}) is equivalent to (\ref{19}). \\
(\ref{3}) follows from (\ref{18}): $2s+1-\rho- \epsilon > 0$. \\
(\ref{4}): $2s+1-\rho-\epsilon+1-r \ge s+1-\frac{1}{p}-2\epsilon > \frac{1}{2} - 
\frac{1}{2p} \ge 0$ by (\ref{19}) and (\ref{18}). \\
(\ref{6}): $2s+\sigma > -1+\frac{2}{p} \ge 0$ by (\ref{18}) and $\sigma > 1/p$. 
\\
(\ref{7}) is weaker than (\ref{8}), which is equivalent to (\ref{19}). \\
(\ref{9}),(\ref{10}),(\ref{11}) are fulfilled, because $1>\sigma > \frac{1}{p}$ 
and $\rho < 1$. \\
(\ref{13}) is satisfied for $\sigma < 1$. \\
(\ref{14}) and (\ref{15}) are equivalent to (\ref{20}). \\
(\ref{16}),(\ref{17}): $\pm s+1-\sigma+r \ge \pm s+1-\sigma+|s| \ge 1-\sigma \ge 
0$ by (\ref{20}). \\
It remains to fulfill 
(\ref{5}): $ r < 2s+\sigma+1-\frac{1}{p}$ and 
(\ref{12}): $r > \sigma - 1 + \frac{1}{p} + \epsilon$ . \\
These conditions can be fulfilled with a suitable $1/p < \sigma < 1$ , provided 
$r<2s+2-\frac{1}{p}$ and $r > \frac{2}{p}-1$ . The last condition is (\ref{21}), 
and the first one is weaker than (\ref{19}) under assumption (\ref{18}):\\
$ 2s+2-\frac{1}{p} > 1+s \, \Longleftrightarrow \, s > -1+\frac{1}{p}$ . \\
This holds, if (\ref{18}) is fulfilled.

Thus we have proven the following
\begin{theorem}
\label{Theorem}
Assume $ 1 < p \le 2$ and (\ref{18}) - (\ref{21}) . The Cauchy problem for the 
Dirac -- Klein -- Gordon system 
(\ref{0.1}),(\ref{0.2}),(\ref{0.3}) with data
$$\psi_0 \in \widehat{H^{s,p}} \,  , \, \phi_0 \in \widehat{H^{r,p}} \, , \, 
\phi_1 \in 
\widehat{H^{r-1,p}} $$
is locally well-posed, i.e. there exists a unique local solution
$$\psi = \psi_+ + \psi_- \quad
{\mbox with} \quad \psi_{\pm} \in X^{s,\sigma}_{\pm p}[0,T] \, ,$$ 
and
$$\phi = \frac{1}{2}(\phi_+ + \phi_-)\, , 
\, \phi_t = \frac{1}{2i} A^{-\frac{1}{2}}(\phi_+ - \phi_-) \quad
{\mbox with} \quad \phi_{\pm} \in Y^{r,\rho}_{\pm p}[0,T] \, ,$$  
where $A=-\frac{\partial^2}{\partial x^2} +1$ and $\rho,\sigma = \frac{1}{p}+$ . 
Here 
$T=T(\|\psi_0\|_{\widehat{H^{s,p}}},\|\phi_0\|_{\widehat{H^{r,p}}} 
,\|\phi_1\|_{\widehat{H^{r-1,p}}})$. This solution satisfies 
$$\psi \in C^0([0,T],\widehat{H^{s,p}}) \, , \, \phi \in 
C^0([0,T],\widehat{H^{r,p}}) 
\, , \, \phi_t \in C^0([0,T],\widehat{H^{r-1,p}}) \, ,$$
and the mapping data upon solution is locally Lipschitz continuous. 
\end{theorem}
{\bf Remark:} From the scaling point of view the spaces $\widehat{H^{s,p}}$ 
behave like the Sobolev spaces $H^{s,p}$ and like $H^{\sigma,2}$ , where $\sigma 
= s + \frac{1}{2} - \frac{1}{p} $ . Similarly $\widehat{H^{r,p}}$ behaves like 
$H^{\lambda,2}$ , where $\lambda = r + \frac{1}{2} - \frac{1}{p} $ . This has 
the following consequences in view of our assumptions (\ref{18}) - (\ref{21}).

Minimizing $\sigma$ requires to take  $s=-\frac{1}{2}+\frac{1}{2p}+$ , so that 
$\sigma_{min} = - \frac{1}{2p}+$ , which is optimal for $p=1+$ , namely 
$\sigma_{min} = - \frac{1}{2}+$ . The corresponding $\lambda$ in the case 
$s=-\frac{1}{2}+\frac{1}{2p}+$ is limited by the conditions (\ref{20}) and 
(\ref{21}), which require $ r > \frac{2}{p} -1$ for $1 < p \le 5/3$ , and $ r > 
\frac{1}{2} - \frac{1}{2p}-$ for $ 5/3 < p \le 2$ . Thus $\lambda_{min} = 
\frac{1}{p} - \frac{1}{2}+$ for $1<p \le 5/3$ , and $ \lambda_{min} = 
1-\frac{3}{2p}-$ for $ 5/3 < p \le 2$, especially for $ p=1+$ : $ \lambda_{min} 
= \frac{1}{2}+$ .

The minimal $\lambda$ , however, is obtained for $s=0$ and $p=2$, namely $r=0+$ 
, $p=2$ , thus $\lambda = 0+$ .

Summarizing, we are able to improve the result from the scaling point of view 
for the Dirac part by leaving the case $p=2$, whereas for the wave part no 
improvement can be given.
\section{Appendix}
We use the notation $ \|f\|_{\widehat{L^p_x}(\widehat{L^q_t})} := 
\|\tilde{f}\|_{L^{p'}_{\xi}(L^{q'}_{\tau})}$ , where $\tilde{f}$ is the Fourier 
transform with respect to space and time. First we prove an embedding theorem.
\begin{prop}
\label{Proposition}
Let $ 1 < r,w_1,w_2,w_3 < \infty$. Then the following embeddings hold:
\begin{eqnarray}
\label{2.1}
X^{\frac{1}{r}+,0}_r & \subset & \widehat{L^{\infty}_x}(\widehat{L^r_t}) \\
\label{2.2}
X^{0,\frac{1}{r}+}_r & \subset & \widehat{L^r_x}(\widehat{L^{\infty}_t}) \\
\label{2.3}
X^{\frac{1}{r}+,\frac{1}{r}+}_r & \subset & 
\widehat{L^{\infty}_x}(\widehat{L^{\infty}_t}) \\
\label{2.4}
X^{\frac{1}{w_1}+,\frac{1}{w_2}+}_{w_3} & \subset & 
\widehat{L^{v_1}_x}(\widehat{L^{v_2}_t}) \, ,
\end{eqnarray}
where $ 1/v_1 = 1/w_3 - 1/w_1 $ , $  1/v_2 = 1/w_3 - 1/w_2 $ .
\end{prop}  
{\bf Proof:}
(\ref{2.1}): H\"older's inequality gives
$$ \|u\|_{\widehat{L^{\infty}_x}(\widehat{L^r_t})} = 
\int\left((\int|\tilde{u}(\xi,\tau)|^{r'} d\tau)^{1/r'} \langle \xi 
\rangle^{1/r+}\right) \langle \xi \rangle^{-1/r-} d\xi \le c 
\|u\|_{X^{1/r+,0}_r} \, . $$
(\ref{2.2}): Similarly
\begin{eqnarray*}
\|u\|_{\widehat{L^r_x}(\widehat{L^{\infty}_t})} & \hspace{-0.8em}= 
&\hspace{-0.8em} \left(\int(\int|\tilde{u}(\xi,\tau)| \langle \tau + \phi(\xi) 
\rangle^{1/r+} \langle \tau + \phi(\xi) \rangle^{-1/r-} d\tau)^{r'} d\xi 
\right)^{1/r'}  \\
&\hspace{-0.8em} \le & \hspace{-0.8em}\left( (\int|\tilde{u}(\xi,\tau)|^{r'} 
\langle \tau + \phi(\xi) \rangle^{r'(\frac{1}{r}+)} d\tau)(\int \langle \tau + 
\phi(\xi) \rangle^{(-\frac{1}{r}-)r} d\tau)^{\frac{r'}{r}} d\xi \right)^{1/r'}  
\\
&\hspace{-0.8em} \le & \hspace{-0.8em}c \|u\|_{X^{0,1/r+}_r} \, .  
\end{eqnarray*} 
(\ref{2.3}): Using H\"older's inequality twice we get
\begin{eqnarray*}
\|u\|_{\widehat{L^{\infty}_t}(\widehat{L^{\infty}_x})} & = & \int \int 
|\tilde{u}(\xi,\tau) \langle \tau + \phi(\xi) \rangle^{1/r+} \langle \tau + 
\phi(\xi) \rangle^{-1/r-} d\tau d\xi \\
& \le & c \int(\int |\tilde{u}(\xi,\tau) \langle \tau + \phi(\xi) 
\rangle^{1/r+}|^{r'} d\tau)^{1/r'} \langle \xi \rangle^{1/r+} \langle \xi 
\rangle^{-1/r-} d\xi \\
& \le & c \|u\|_{X^{1/r+,1/r+}_r}
\end{eqnarray*}
(\ref{2.4}): In a first step we interpolate between (\ref{2.1}) and (\ref{2.3}) 
and get
\begin{equation}
\label{2.5}
X^{1/r+,1/r_1 +}_r \subset \widehat{L^{\infty}_x}(\widehat{L^{r_2}_t}) \, ,
\end{equation}
where $1/r_2 = 1/r - 1/r_1$ . Similarly, interpolation between (\ref{2.2}) and 
(\ref{2.3}) gives
\begin{equation}
\label{2.6}
X^{1/r_1 +,1/r +}_r \subset \widehat{L^{r_2}_x}(\widehat{L^{\infty}_t}) \, ,
\end{equation}
where $1/r_2 = 1/r - 1/r_1$ . Finally, interpolating between (\ref{2.5}) and 
(\ref{2.6}) gives the desired result.

We now use this proposition to prove a product law for $X^{s,b}_{r,\phi}$ - 
spaces, belonging to the phase function $\phi$.
\begin{prop}
\label{Proposition2.1}
Let $ 1 < p,q,r < \infty $ , $ \alpha, \beta, \gamma \ge 0 $ . Assume $ a+b+c > 
1/p + 1/q + 1/r' - 1 $ , $ a+b \ge 0$ , $a+c \ge 0$ , $b+c \ge 0$ and $\alpha + 
\beta + \gamma > 1/p + 1/q + 1/r' -1. $ Then the following estimate holds:
$$ \|uv\|_{X^{-c,-\gamma}_{r,\phi}} \le c \|u\|_{X^{a,\alpha}_{p,\phi_1}} 
\|v\|_{X^{b,\beta}_{q,\phi_2}} \, . $$
\end{prop}
{\bf Proof:} Using $ (X^{-c,-\gamma}_{r,\phi})^* = X^{c,\gamma}_{r',\phi}$ we 
have to show
$$ \left| \int\int_* \frac{\tilde{v}_1(\xi_1,\tau_1)}{\langle \xi_1 \rangle^a 
\langle \sigma_1 \rangle^{\alpha}} \frac{\tilde{v}_2(\xi_2,\tau_2)}{\langle 
\xi_2 \rangle^b \langle \sigma_2 \rangle^{\beta}} 
\frac{\tilde{\varphi}(\xi,\tau)}{\langle \xi \rangle^c \langle \sigma 
\rangle^{\gamma}} d\xi_1 d\xi_2 d\tau_1 d\tau_2 \right| \le c 
\|v_1\|_{\widehat{L^p}} \|v_2\|_{\widehat{L^q}}\|\varphi\|_{\widehat{L^{r'}}} \, 
, $$
where $\sigma_i = \tau_i + \phi_i(\xi_i)$ , $\sigma = \tau + \phi(\xi)$ . \\
Assume $a,b,c \ge 0$ first. Applying Young's and H\"older's inequality we get a 
bound
$$c \|J^{-a} \Lambda_1^{-\alpha} v_1 
\|_{\widehat{L^{p_1}_x}(\widehat{L^{p_2}_t})} \|J^{-b} \Lambda_2^{-\beta} v_2 
\|_{\widehat{L^{q_1}_x}(\widehat{L^{q_2}_t})} \|J^{-c} \Lambda^{-\gamma} \varphi 
\|_{\widehat{L^{r_1}_x}(\widehat{L^{r_2}_t})} \, $$
where $1/p_1 + 1/q_1 + 1/r_1 = 1 = 1/p_2 + 1/q_2 + 1/r_2$ . Here $J$ and 
$\Lambda_1 , \Lambda_2 , \Lambda$ are the Fourier multipliers with symbols 
$\langle \xi \rangle$ and $\langle \sigma_1 \rangle , \langle \sigma_2 \rangle , 
\langle \sigma \rangle$ , respectively. In order to get the desired bound we 
need the embeddings:
$$ X^{a,\alpha}_{p,\phi_1} \subset \widehat{L^{p_1}_x}(\widehat{L^{p_2}_t}) \, , 
\,
X^{b,\beta}_{q,\phi_2} \subset \widehat{L^{q_1}_x}(\widehat{L^{q_2}_t}) \, , \, 
X^{c,\gamma}_{r',\phi} \subset \widehat{L^{r_1}_x}(\widehat{L^{r_2}_t}) \, . $$
These embeddings hold by Proposition \ref{Proposition}, if the following 
conditions are satisfied:
\begin{eqnarray*}
 1/p_1 > 1/p - a \, , \, 1/p_2 > 1/p - \alpha \, , \, 1/q_1 > 1/q - b \\ 
1/q_2 > 1/q - \beta \, , \, 1/r_1 > 1/r' - c \, , \, 1/r_2 > 1/r' - \gamma \, . 
\end{eqnarray*}
These inequalities can obviously be fulfilled if
$1=1/p_1 + 1/q_1 + 1/r_1 > 1/p + 1/q + 1/r' -(a+b+c)  \Leftrightarrow  a+b+c > 
1/p+1/q+1/r'-1$ and
$1=1/p_2 + 1/q_2 + 1/r_2 > 1/p + 1/q + 1/r' -(\alpha +\beta +\gamma)  
\Leftrightarrow  \alpha + \beta + \gamma > 1/p+1/q+1/r'-1$ . \\
Consider now the case, where not all the numbers $a,b,c$ are nonnegative. In 
view of our assumptions only one number is negative, $a<0$ , say. Then $\langle 
\xi_1 \rangle^{-a} \le c(\langle \xi_2 \rangle^{-a} + \langle \xi \rangle^{-a})$ 
. Thus we only have to consider the integral with $(a,b,c)$ replaced by 
$(0,a+b,c)$ or $(0,b,a+c)$. All these numbers are nonnegative and their sum also 
fulfills our assumptions. This proves the general case.

Finally, we prove a bilinear estimate for free waves.
\begin{prop}
\label{Proposition2.2}
If $u$ and $v$ are solutions of the linear problems
\begin{eqnarray*}
(\frac{\partial}{\partial t} + \frac{\partial}{\partial x})u = 0 & , & u(x,0) = 
f(x) \, , \\
(\frac{\partial}{\partial t} - \frac{\partial}{\partial x})v = 0 & , & v(x,0) = 
g(x) \, ,
\end{eqnarray*}
where $f,g \in \widehat{L^p}$ , $ 1 < p < \infty $ , the following estimate 
holds:
$$ \|uv\|_{\widehat{L^p}} \le c \|f\|_{\widehat{L^p}} \|g\|_{\widehat{L^p}} \, . 
$$
\end{prop}
{\bf Proof:} We have $u(x,t) = f(x-t)$ and $v(x,t) = g(x+t)$ and thus 
$\tilde{u}(\xi,\tau) = \delta(\tau + \xi) \widehat{f}(\xi)$ and 
$\tilde{v}(\xi,\tau) = \delta(\tau - \xi) \widehat{g}(\xi)$ . This implies
\begin{eqnarray*}
{\cal F}(uv)(\xi,\tau) & = & \int\int \delta(\lambda + \eta) \widehat{f}(\eta) 
\delta(\tau - \lambda -(\xi-\eta)) \widehat{g}(\xi - \eta) d\lambda d\eta \\
& = & \int \int \delta(\tau +2\eta - \xi) \widehat{f}(\eta) \widehat{g}(\xi - 
\eta) d\eta \\
& = & \widehat{f}\left(\frac{\xi - \tau}{2}\right) \widehat{g}\left(\frac{\xi + 
\tau}{2}\right) \, .
\end{eqnarray*}
Consequently
$$ \|uv\|_{\widehat{L^p}} = \left(\int\int\left|\widehat{f}\left(\frac{\xi - 
\tau}{2}\right) \widehat{g}\left(\frac{\xi + \tau}{2}\right)\right|^{p'} d\xi 
d\tau \right)^{1/p'} \le c \|f\|_{\widehat{L^p}}  \|g\|_{\widehat{L^p}} \, . $$ 

This proposition implies
\begin{Cor}
\label{Corollary2.1}
Let $1<p<\infty$ , $\sigma > 1/p$ . Then
$$ \|uv\|_{\widehat{L^p}} \le c \|u\|_{X^{0,\sigma}_{+p}} 
\|v\|_{X^{0,\sigma}_{-p}} \, . $$
\end{Cor}
{\bf Proof:} Recall that (with $U_{\pm}(t)$ denoting the evolution operator of 
the equation $(\frac{\partial}{\partial t} \pm \frac{\partial}{\partial 
x})u=0$):
$$ \|u\|_{X^{0,\sigma}_{\pm p}} =  (\int \langle \tau \rangle^{\sigma p'} 
\|{\cal F}(U_{\pm}(-\cdot)u)(\xi,\tau)\|_{L^{p'}_{\xi}}^{p'} d\tau)^{1/p'} \, . 
$$  
We use
$$ u(t) = c \int e^{it\tau} U_+(t) h(\tau) d\tau \quad {\mbox with} \quad 
h:={\cal F}_t U_+(- \cdot)u \,$$ 
and
$$ v(t) = c \int e^{it\rho} U_-(t) l(\rho) d\rho \quad {\mbox with} \quad 
l:={\cal F}_t U_-(- \cdot)u \, .$$
Thus
$$ (uv)(t) = c \int \int e^{it\tau} U_+(t)h(\tau) e^{it\rho} U_-(t)l(\rho) d\tau 
d\rho \, . $$
By Minkowski's inequality
$$ \|uv\|_{\widehat{L_{xt}^p}} \le c \int \int \|e^{it(\tau + \rho)} 
U_+(t)h(\tau) U_-(t)l(\rho)\|_{\widehat{L_{xt}^p}} d\tau d\rho \, . $$
But now for fixed $\tau$ and $\rho$ we have
\begin{eqnarray*}
\|e^{it(\tau + \rho)} w(x,t)\|_{\widehat{L_{xt}^p}} & = & \| {\cal F}(e^{it(\tau 
+ \rho)}w(x,t))\|_{L^{p'}_{\xi \tau'}} \\ & = & c\| \int \int e^{-i(x\xi 
+t(\tau' - \tau - \rho))}w(x,t) dx dt \|_{L^{p'}_{\xi \tau'}} \\
& = & c\|\tilde{w}(\xi,\tau' - \tau - \rho)\|_{L^{p'}_{\xi \tau'}} = 
c\|\tilde{w}(\xi,\tau')\|_{L^{p'}_{\xi \tau'}} = c\|w\|_{\widehat{L^p_{xt}}} \, 
,
\end{eqnarray*}
so that with $w(t)=U_+(t)h(\tau)U_-(t)l(\rho)$ we get by Proposition 
\ref{Proposition2.2}:
\begin{eqnarray*}
\|uv\|_{\widehat{L^p_{xt}}} & \le & c \int\int 
\|U_+(t)h(\tau)U_-(t)l(\rho)\|_{\widehat{L^p_{xt}}} d\tau d\rho \\
& \le & c \int \|h(\tau)\|_{\widehat{L^p_x}} d\tau \int 
\|l(\rho)\|_{\widehat{L^p_x}} d\rho \\
& = & c \int\langle \tau \rangle^{-\sigma} \langle \tau \rangle^{\sigma} \|{\cal 
F}(U_+(- \cdot)u)\|_{L^{p'}_{\xi}} d\tau \int\langle \rho \rangle^{-\sigma} 
\langle \rho \rangle^{\sigma} \|{\cal F}(U_-(- \cdot)v)\|_{L^{p'}_{\xi}} d\rho 
\\
& \le & c (\int \langle \tau \rangle^{-\sigma p} d\tau)^{1/p} (\int \langle \tau 
\rangle^{\sigma p'} \|{\cal F}(U_+(- \cdot)u)\|_{L^{p'}_{\xi}}^{p'} 
d\tau)^{1/p'} \cdot \\
& & \quad \cdot (\int \langle \rho \rangle^{-\sigma p} d\rho)^{1/p} (\int 
\langle \rho \rangle^{\sigma p'} \|{\cal F}(U_-(- \cdot)v)\|_{L^{p'}_{\xi}}^{p'} 
d\rho)^{1/p'} \\
& \le & c \|u\|_{X^{0,\sigma}_{+p}} \|v\|_{X^{0,\sigma}_{-p}} \, .
\end{eqnarray*}
We also used the following elementary algebraic inequality for real numbers.
\begin{lemma}
\label{Lemma2.1}
If $\xi_1,\xi_2,\tau_1,\tau_2 \in {\bf R}$ and $\xi = \xi_1 + \xi_2$ , $ \tau = 
\tau_1 + \tau_2$ , the following estimate holds:
$$ \min(|\xi_1|,|\xi_2|) \le 1/2 (|\sigma_ {\pm}| + |\sigma_1^+| + |\sigma_2^-|) 
\, , $$ 
where
$$ \sigma_{\pm} := \tau \pm |\xi| \, , \, \sigma_1^+ := \tau_1 + \xi_1 \, , \, 
\sigma_2^- := \tau_2 - \xi_2 \, . $$
\end{lemma}
{\bf Proof:} We have
$$
\sigma_{\pm}  =  \tau \pm |\xi|  =  \sigma_1^+ + \sigma_2^- \pm |\xi_1 + \xi_2| 
- \xi_1 + \xi_2 \, . $$
Now, if $\xi_1 + \xi_2 {\ge \atop \le} 0$ , we have $\pm |\xi_1 + \xi_2| - \xi_1 
+ \xi_2 = 2\xi_2$ , whereas, if $\xi_1 + \xi_2 { \le \atop \ge} 0$ we have $\pm 
|\xi_1 + \xi_2| - \xi_1 + \xi_2 = -2\xi_1$ , so that the claimed inequality 
follows.


\begin{thebibliography}{A1}
\bibitem[AFS]{AFS} P. d'Ancona, D. Foschi, and S. Selberg: {\sl Null structure 
and 
almost optimal local regularity for the Dirac -- Klein -- Gordon system}. 
arXiv: math. AP/0509545 , to appear in Journal of the EMS
\bibitem[B]{B} N. Bournaveas: {\sl A new proof of global existence for the Dirac 
Klein--Gordon equations in one space dimension}. J. Funct. Analysis 173 (2000), 
203-213
\bibitem[BG]{BG} N. Bournaveas and D. Gibbeson: {\sl Low regularity global 
solutions of the Dirac -- Klein -- Gordon equations in one space dimension}. 
Diff. Int. Equations 19 (2006), 211-222
\bibitem[CVV]{CVV} T. Cazenave, L. Vega and M.C. Vilela: {\sl A note on the 
nonlinear Schr\"odinger equation in weak $L^p$ spaces}. Comm. Contemp. Math. 3 
(2001), 153-162 
\bibitem[C]{C} J. M. Chadam: {\sl Global solutions of the Cauchy problem for the 
(classical) coupled Maxwell -- Dirac equations in one space dimension}. J. 
Funct. Analysis 13 (1973), 173-184
\bibitem[CG]{CG} J. M. Chadam and R. T. Glassey: {\sl On certain global 
solutions 
of the Cauchy problem for the (classical) coupled Klein -- Gordon -- Dirac 
equations in one and three space dimensions}. Arch. Rat. Mech. Anal. 54 (1974), 
223-237
\bibitem[F]{F} Yung-Fu Fang: {\sl A direct proof of global existence for the 
Dirac -- Klein -- Gordon equations in one space dimension}. Taiwanese J. Math. 8 
(2004), 33-41
\bibitem[F1]{F1} Yung-Fu Fang: {\sl Low regularity solutions for Dirac -- Klein 
-- Gordon equations in one space dimension}. Electr. J. Diff. Equ. 
2004(2004), no. 102, 1-19
\bibitem[G1]{G1} A. Gr\"unrock: {\sl An improved local well-posedness result 
for the modified KdV equation}. Int. Math. Res. Not. 2004, no. 61, 3287-3308
\bibitem[G2]{G2} A. Gr\"unrock: {\sl Bi- and trilinear Schr\"odinger estimates 
in one space dimension with applications to cubic NLS and DNLS}. Int. Math. 
Res. Not. 2005, no. 41, 2525-2558
\bibitem[M]{M} S. Machihara: {\sl The Cauchy problem for the 1d 
Dirac-Klein-Gordon equation}. to appear in NoDEA
\bibitem[P]{P} H. Pecher: {\sl Low regularity well-posedness for the 
one-dimensional Dirac-Klein-Gordon system}. Electr. J. Diff. Equ. 2006(2006), 
no. 150, 1-13
\bibitem[ST]{ST} S. Selberg and A. Tesfahun: {\sl Low regularity well-posedness 
for the one dimensional Dirac-Klein-Gordon system}. arXiv: math. AP/0611718
\bibitem[VV]{VV} A. Vargas and L. Vega: {\sl Global wellposedness for 1D 
non-linear Schr\"odinger equation for data with infinite $L^2$ norm}. J. Math. 
Pures Appl. 80 (2001), 1029-1044
\end{thebibliography}
\end{document}